\documentclass[11pt]{article}
\usepackage{amsmath}
\usepackage{amssymb}
\usepackage{txfonts}
\usepackage{mathptmx}     
\usepackage{latexsym}
\usepackage{relsize}
\newtheorem{teo}{Theorem}
\newtheorem{lem}{Lemma}
\newtheorem{rem}{Remark}
 
 \newtheorem{coro}{Corollary}
 \newtheorem{conje}{Conjecture}
 \newcommand{\adots}
 {\mathinner{\mkern2mu\raise1pt\hbox{.}\mkern2mu\raise4pt\hbox{.}\mkern2mu\raise7pt\hbox{.}\mkern1mu}}
\title{The bispectral problem and polynomial solutions}

\author{D. Barrios Rolan\'{\i}a\\
Dpto. Matem\'atica Aplicada a la Ingenier\'{\i}a Industrial\\
ETSI Industriales\\
Universidad Polit\'ecnica de Madrid\\
Spain\\}

\date{}
\begin{document}

\maketitle

\begin{abstract}
Some properties and relations satisfied by the polynomial solutions of the bispectral problem are studied. Given a differential operator, under certain restrictions its polynomial eigenfunctions are explicitly obtained, as well as the corresponding eigenvalues. Also the inverse problem is approached, giving conditions for the existence of a differential operator from its eigenpolynomials and eigenvalues. 
\end{abstract}

\section{Introduction}

We consider the ordinary differential operator
\begin{equation}
L\equiv                
\sum_{i=0}^Na_i(x)\partial_x^i
\label{1}
\end{equation}
where 
$a_i(x)
$
are polynomials,
$
\deg (a_i)\leq i,
$ and 
$
\partial^i_x=\frac{d^i}{dx^i},\, i=1,\ldots,N
$. We also consider a sequence 
$
\{\lambda_n\}\subset \mathbb{C}
$
of its eigenvalues and the corresponding sequence of eigenfunctions
$
\{P_n\},
$
which we assume that are monic polynomials with $\deg(P_n)=n$ for each $n\in \mathbb{N}$. This is,
\begin{equation}        
\sum_{i=1}^Na_i(x)\partial_x^iP_n(x)=\lambda_nP_n(x)\,,\quad \forall n\in \mathbb{N}.
\label{11}
\end{equation}
The polynomials $P_n$ verifying \eqref{11} are called {\em eigenpolynomials} in this work. The Sturm-Liouville polynomial problem deals with the relations between the operator $L$ and its eigenvalues and eigenpolynomials. We are interested in families of such eigenpolynomials that are at the same time eigenfunctions of a certain difference operator $D$. More precisely, 
\begin{equation}
\label{diferencias}
D
\left(
\begin{array}{c}
P_0(x)\\
P_1(x)\\
\vdots \\
\end{array}
\right)=
x\left(
\begin{array}{c}
P_0(x)\\
P_1(x)\\
\vdots \\
\end{array}
\right),
\end{equation}
where
$$
D=
\left(
\begin{array}{rcccccccccc}
\alpha_{0,0}&1&0&\cdots  \\ 
 \vdots&\ddots&\ddots&\ddots  \\
 \alpha_{p,0}&\cdots & \alpha_{p,p}&1&0&\cdots \\ 
0&\alpha_{p+1,1}&\cdots&\alpha_{p+1,p+1}&1&0& \cdots\\ 
 & 0&\ddots&  \ddots&  \ddots& \ddots&\\
&&0&\alpha_{n,n-p}&\cdots&\alpha_{n,n}&1&0& \cdots\\ 
 &&&\ddots &\ddots&  \ddots& \ddots& \ddots& \\
 \end{array}
\right).
$$
This is, the polynomials
$
\{P_n\},\,n\in\mathbb N,
$ 
verify a $(p+2)$-term recurrence relation
\begin{equation}
\label{otro9}
\displaystyle \sum_{k=n-p}^{n-1}\alpha_{n,k}P_k(x)+(\alpha_{n,n}-x)P_n(x)+P_{n+1} = 0\,,\quad n=0,1,\ldots ,
\end{equation}
with initial conditions 
$$
P_0=1,\quad P_{-1}=\cdots =P_{-p}=0.
$$
The difference operator are given by 
$$
D(n)P_n=\sum_{k=n-p}^{n}\alpha_{n,k}P_k(x)+P_{n+1}
$$
and we have 
\begin{equation}
\label{discreto}
D(n)P_n=xP_n.
\end{equation}
If the polynomials
$
P_n
$ 
verify \eqref{11} and \eqref{discreto}, we say that the sequence
$
\{P_n\},\,n\in\mathbb N,
$ 
is a solution for the bispectral problem defined by $L$ and $D$.

S. Bochner \cite{Bothner} studied the above problem when the differential operators has order $N=2$, determining such polynomial solutions. He solved completly the problem, and his classification defined the well-known nowaday families of classical orthogonal polynomials. Some years after, Krall \cite{Krall1}, \cite{Krall} deals the differential operator of order $N=4$, giving a classification containing seven families, including three new families of polynomial eigenfunctions that do not reduce to operators of order 2. 

Since the celebrated paper of Bochner \cite{Bothner} and the relevant contributions of Krall, much has been written about the bispectral problem \cite{Ben}, \cite{Shapiro}, \cite{Zhedanov}. We highlight the paper \cite{horo} for its excellent review of the state of the art, as well as for its contributions to operators of order $N = 3$. Other extensions of these works have been attempted, even for operators in differences with complex coefficients \cite{Grunbaum1}, \cite{Grunbaum2}. However, the hardness of the problem has so far not made it possible to obtain conclusions as relevant as those already known for $N = 2$ and $N=4$. 

The goal of this work is to shed some light on this problem providing some relationships between the differential operator \eqref{1} and its eigenpolynomials. We hope that these contributions open new possibilities for solving the problem. In Section \ref{seccion2} the main result is Theorem \ref{lema1}, where conditions on the coefficients of the eigenpolynomials of a difference operator are stated. Moreover other auxiliary results used in the rest of the paper are showed. In Section \ref{seccion3} we assume that the differential operator $L$ is given and then its eigenvalues and eigenpolynomials are constructed. Finally, in Section \ref{seccion4} the inverse problem is approached, setting the differential operator $L$ from its eigenvalues and eigenpolynomials. 

\section{Eigenpolynomials and eigenvalues\label{seccion2}}
In \eqref{1} we assume $a_0\equiv 0$ because, on the contrary, we substitute $a_0$ by $a_0-\lambda_0$. For the same reason, we assume $\lambda_0=0$. 

Writing
\begin{equation}
a_m(x)=\sum_{i=0}^m a_{m,i}x^i,\quad m=1,\ldots, N,
\label{1111}
\end{equation}
and
\begin{equation}
P_n(x)=\sum_{i=0}^n b_{n,i}x^i,\quad b_{n,n}=1, \quad n\in \mathbb{N},
\label{111}
\end{equation}
in the sequel we analyze the relation between the coefficients of $P_n(x)$ and those of $a_m(x)$ when \eqref{11} is verified. However, previously some auxiliary results of a more general character are showed.

\begin{lem}
\label{combinatorio}
For all $r\in \mathbb{N}$ we have
$$
\sum_{i=k}^{r+1}(-1)^{i-k}{i \choose k}{r+1 \choose i}=0\,,\quad k=0,1,\ldots , r\,.$$
\end{lem}
(The proof of this lemma is an exercise of elemental Combinatorial Analysis.)

\begin{lem}
\label{lema3}
Given two sequences $\{a_{n,i}\},\,n=0,1,\ldots ,\,i=0,1,\ldots, n\,,$ and $\{\delta_n^{(k)}\},\,n=0,1,\ldots ,\,k=0,1,\ldots, n\,,$ the following connections are equivalent.
\begin{itemize}
\item[1.-] For each  $n=0,1,\ldots ,$
\begin{equation}
\label{A}
\delta_n^{(k)}=\sum_{i=k}^{n} {n \choose i} i!a_{i,i-k}\quad \text{ for } k=0,1,\ldots, n\,.
\end{equation}

\item[2.-] For each  $n=0,1,\ldots ,$
\begin{equation}
\label{B}
n!a_{n,n-k}=\sum_{i=k}^{n} {n \choose i}(-1)^{n-i}\delta_i^{(k)}\quad \text{ for } k=0,1,\ldots, n\,.
\end{equation}

\end{itemize}

\end{lem}
\noindent
\underline{Proof.}- We proceed by induction on $n$. 

For $n=0$ we have $a_{0,0}=\delta_0^{(0)}$ in \eqref{A} and \eqref{B}.

Suppose \eqref{A} holds for $n=0,1,\ldots$ and assume that \eqref{B} is verified for $n=0,1,\ldots, m$, with $m\in \mathbb{N}$. Taking $n=m+1$ in \eqref {A}, we have
$$
\delta_{m+1}^{(k)}=\sum_{i=k}^{m+1} {m+1 \choose i} i!a_{i,i-k}=(m+1)!a_{m+1,m+k-1}+\sum_{i=k}^{m} {m+1 \choose i} i!a_{i,i-k}
$$
for 
$k=0,1,\ldots, m+1$. Therefore,
\begin{eqnarray}
\label{7tilde}
(m+1)!a_{m+1,m+k-1}&= & \delta_{m+1}^{(k)}-\sum_{i=k}^{m} {m+1 \choose i} \left[
\sum_{s=k}^{i} {i \choose s}(-1)^{i-s}\delta_s^{(k)}\right]\nonumber \\
& = & \delta_{m+1}^{(k)}-\sum_{s=k}^{m} \left[ 
\sum_{i=s}^{m} {m+1 \choose i} {i \choose s}(-1)^{i-s}\right] \delta_s^{(k)}\,.
\end{eqnarray}
Since Lemma \ref{combinatorio},
$$
\sum_{i=s}^{m} {m+1 \choose i} {i \choose s}(-1)^{i-s}=-{m+1 \choose s} (-1)^{m+1-s}\,.
$$
Then, from \eqref{7tilde},
$$
(m+1)!a_{m+1,m+k-1}= \delta_{m+1}^{(k)}+\sum_{s=k}^{m} 
{m+1 \choose s} (-1)^{m+1-s}\delta_s^{(k)}\,,
$$
which drives to \eqref{B} for $n=m+1$. 

Conversely, assuming that \eqref{B} is verified for $n=0,1,\ldots $ and \eqref{A} holds for $n=0,1,\ldots ,m$, and taking $n=m+1$ in \eqref{B}, 
\begin{eqnarray*}
\delta_{m+1}^{(k)}&= &(m+1)!a_{m+1,m+k-1} -\sum_{i=k}^{m} {m+1 \choose i}(-1)^{m+1-i}\delta_i^{(k)}\\
& = &(m+1)!a_{m+1,m+k-1} -\sum_{i=k}^{m} {m+1 \choose i}(-1)^{m+1-i} \left[ 
\sum_{s=k}^{i} {i \choose s} s!a_{s,s-k}\right]\\
& = &(m+1)!a_{m+1,m+k-1} -\sum_{s=k}^{m} \left[ 
\sum_{i=s}^{m} {m+1 \choose i} {i \choose s}(-1)^{m+1-i} \right]s!a_{s,s-k}\\
& = &(m+1)!a_{m+1,m+k-1} +\sum_{s=k}^{m} (-1)^{m-s} \left[ 
\sum_{i=s}^{m} {m+1 \choose i} {i \choose s}(-1)^{i-s} \right]s!a_{s,s-k}\,.
\end{eqnarray*}
Using again Lemma \ref{combinatorio},
$$
\sum_{i=s}^{m} {m+1 \choose i} {i \choose s}(-1)^{i-s} =-{m+1 \choose s}(-1)^{m+1-s} \,.
$$
Hence,
$$
\delta_{m+1}^{(k)}=(m+1)!a_{m+1,m+k-1} +\sum_{s=k}^{m}{m+1 \choose s} s!a_{s,s-k}\,,
$$
which gives \eqref{A} for $n=m+1$. \hfill $\square$

We underline that the above auxiliary results are independent on the fact that the sequences 
$
\{\delta_n^{(k)}\}
$
and 
$
\{a_{n,n-k}\}
$
are related with the differential operator $L$. However, in the sequel we deal with the sequences 
$
\{a_{i,j}\}
$
and
$
\{b_{i,j}\}
$
that define the polynomials \eqref{1111}--\eqref{111} corresponding to the coefficients of the differential operator $L$ in \eqref{1}  and its eigenpolynomials $P_n$ verifying \eqref{11}, respectively. The following theorem gives the key to understanding the relationship between these sequences of coefficients. 

\begin{teo}
\label{lema1}
Consider the polynomials 
$
a_m(x),\,m=1,\ldots,N,
$
and the sequence 
$
\{P_n(x)\},\, n\in\mathbb N,
$
as in \eqref{1111}--\eqref{111}. Then \eqref{11} is verified if and only if 
\begin{equation}
\sum_{k=0}^N\left[\sum_{i=k}^{m+k} {m+k \choose i} i!a_{i,i-k}\right]b_{n,m+k}=\lambda_nb_{n,m},\quad n\geq 0\,,\quad m=0,1,\ldots, n\,.
\label{3}
\end{equation}
In particular, if \eqref{11} takes place then the eigenvalues are obtained as
\begin{equation}
\lambda_n=\sum_{i=1}^{\min\{n,N\}} 
{n \choose i} i!a_{i,i}.
\label{2}
\end{equation}
\end{teo}

\noindent
\underline{Proof.}- We follow the lines of Bocher in \cite{Bothner}. For each $n\in \mathbb{N}$, since \eqref{111}, 
\begin{eqnarray*}
\partial_x^{i}P_n(x)& = & \sum_{s=1}^{n-i+1}(i+s-1)(i+s-2))\cdots sb_{n,i+s-1}x^{s-1}\\
& = & \sum_{s=1}^{n-i+1}{i+s-1 \choose i}i!b_{n,i+s-1}x^{s-1},\quad i=1,\ldots, N,
\end{eqnarray*}
where $b_{n,j}=0$ if $j>n$. Then, \eqref{11} is equivalent to
\begin{equation}
\sum_{i=1}^N\left(\sum_{k=0}^i a_{i,k}x^k\right) \left(\sum_{s=1}^{n-i+1}{i+s-1 \choose i}i!b_{n,i+s-1}x^{s-1}\right)=\lambda_n\sum_{r=0}^n b_{n,r}x^r.
\label{22}
\end{equation}
Comparing the coefficients of $x^r, r=0,\ldots,$ in both sides of \eqref{22} we arrive to \eqref{3}.

Finally, if $m=n$ in \eqref{3}, then $b_{n,m+k}=0$ for $k\geq 1$ on the left hand side of \eqref{3}  and we obtain \eqref{2}.  \hfill $\square$

\section{The direct problem\label{seccion3}} 

In this section we assume that the differential operador \eqref{1} is given.  Then 
$N\in \mathbb N$
is fixed and $a_m(x)$ verifies \eqref{1111} for each $m=1,\ldots,N$. We take
\begin{equation}
a_m(x)=0,\quad m>N.
\label{aaa}
\end{equation}
Our aim is to find the sequence of eigenvalues 
$\{\lambda_n\}$
and the corresponding sequence of eigenpolynomials
$
\{P_n\} 
$
verifying \eqref{11}. If there exist such sequences, then \eqref{3} holds and \eqref{2} defines the sequence $\{\lambda_n\}$ of eigenvalues. We remark the uniqueness of $\{\lambda_n\}$  for each differential operator $L$. 

Hence, we focus our interest on obtaining the eigenpolynomials. The following result provides a valuable method to obtain the coefficients $b_{n,k},\,k=0,1,\ldots ,n-1,$ of the searched polynomials $P_n(x)$ (see \eqref{111}).

In this section, from the differential operator \eqref{1}, for each $n\in\mathbb N$ we define $\delta_n^{(k)}$ as in \eqref{A}. 

\begin{teo} 
\label{teorema1}
Let $M$ be the semi-infinite upper triangular matrix 
\begin{equation}
M=\left(\begin{array}{ccccccccc}
\delta^{(0)}_{0}&\delta^{(1)}_{1}&\cdots&\cdots&\delta^{(N)}_{N}&0&&&\\
 0&\delta^{(0)}_{1}&\delta^{(1)}_{2}&\cdots&\cdots&\delta^{(N)}_{N+1}&0&&\\
 &0&\delta^{(0)}_{2}&\delta^{(1)}_{3}&\cdots&\cdots&\delta^{(N)}_{N+2}&0&\\
 & &\ddots&\ddots&\ddots&\cdots&&\ddots&\ddots
\end{array}\right)
\end{equation}
and let $M_{n+1}$ be its main section of order $n+1$. Assume 
\begin{equation}
\lambda_n\neq 0, \lambda_1,\,\lambda_2,\ldots, \lambda_{n-1},
\label{autovectores}
\end{equation}
where 
$
\{\lambda_n\}
$
is defined as in \eqref{2}. Then there exists a unique monic polynomial $P_n(x)$ verifying \eqref{11} whose coefficients $b_{n,0},\cdots, b_{n,n-1}$ in \eqref{111} determine an eigenvector of $M_{n+1}$ corresponding to the eigenvalue $\lambda_n$. This is, 
\begin{equation}
\label{otro11}
(M_{n+1}-\lambda_nI_{n+1})b_n=0,
\end{equation}
where $b_n=(b_{n,0},\cdots, b_{n,n-1},1)^T$. 
\end{teo}

\noindent
\underline{Proof}.- Since \eqref{aaa} we have
$
a_{i,i-k}=0
$
for $i>N$. Then, from \eqref{A}, 
$$
\delta_n^{(k)}=\sum_{i=k}^{\min\{n,N\}} {n \choose i} i!a_{i,i-k},\quad k=0,1,\ldots, n\,.
$$
In particular, in \eqref{2} we have 
\begin{equation}
\lambda_n=\delta_n^{(0)},\quad n=0,1,\ldots 
\label{****}
\end{equation}

Due to \eqref{A}, we can rewrite \eqref{3} as
\begin{equation}
\sum_{k=0}^N\delta_{m+k}^{(k)} b_{n,m+k}=\lambda_nb_{n,m},\quad m=0,1,\ldots, n,
\label{7}
\end{equation}
where $b_{n,m+k}=0$ when $m+k>n$. This is, for each fixed $n$,
\begin{equation}
\label{matriz}
(\delta^{(0)}_{m}-\lambda_n)b_{n,m}+\delta^{(1)}_{m+1}b_{n,m+1}+\cdots+ \delta^{(N)}_{m+N}b_{n,m+N}=0\,,\quad m=0,1,\ldots , n,
\end{equation}
assuming $ b_{n,j}=0$ for $j>n$. \eqref{matriz} can be interpreted as an upper triangular linear system, with matrix of coefficients $M_{n+1}$, whose unknowns are $ b_{n,0},\cdots, b_{n,n-1},b_{n,n}$. In matrix notation this is \eqref{otro11}. Hence the coefficients of $b_n$ define an eigenvector of $M_{n+1}$ corresponding to the eigenvalue $\lambda_n$.

Since \eqref{****} and the fact that 
$
\det(M_n-\lambda_n I_n)\neq 0\,,
$
the coefficients $ b_{n,0},\cdots, b_{n,n-1}$ are uniquely determined. 
\hfill $\square$

In the above proof we see that condition \eqref{autovectores} is not neccesary for the existence of eigenpolynomials. However,  under such condition Theorem \ref{teorema1} provides an easy method to obtain the coefficients of these polynomials as coordinates of eigenvectors for a sequence of finite triangular matrices. Henceforth in the sequel we assume that the hypotheses of the Theorem \ref{teorema1} are satisfied. This is, we assume that the eigenvalues 
$\lambda_n,\,n=0,1,\ldots,  $ are all different from each other, which depends on the leading coefficients
$
a_{m,m}
$
in \eqref{1111}. 

\begin{rem}
\label{remark 1}
Because 
$
\delta_n^{(k)}\,,\,k=0,1,\ldots ,n,
$
are polynomials of $n$, from \eqref{7} we see the well-known fact that the coefficients 
$
b_{n,j}
$
are rational functions of $n$.
\end{rem}

We note that \eqref{matriz} is the extension of expresion (11) in \cite{Bothner} to the case of differential operators \eqref{1} of order $N>2$. This provides the explicit expression of each eigenpolynomial, as is showed in the following result.

\begin{teo}
\label{teorema_coeficientes}
Under the above conditions, for each $n\in\mathbb{N}$ we have 
\begin{equation}
\label{coeficientes}
b_{n,i}=\mathlarger\sum_{E_n^{(i)}}\left(\prod_{s=1}^k\frac{\delta^{(i_s)}_{i+i_1+\cdots +i_s}}{\lambda_n-\lambda_{i+i_1+\cdots +i_{s-1}}}\right),\quad i=0,1,\ldots , n-1,
\end{equation}
where the sum is extended to the set
$
E_n^{(i)}=\{(i_1,\ldots,i_k)\,:\, k\in \mathbb N, \,i_1+\cdots +i_k=n-i\}
$
and we understand 
$
i+i_1+\cdots +i_{s-1}=i
$
when 
$s=1 $.
\end{teo}
\noindent
\underline{Proof.}- Taking $i=n-m$ in \eqref{matriz} we see 
$$
\delta^{(0)}_{m}b_{m+i,m}+\delta^{(1)}_{m+1}b_{m+i,m+1}+\cdots+ \delta^{(N)}_{m+N}b_{m+i,m+N}=\lambda_{m+i}b_{m+i,m},
$$
where, as usual, 
$
b_{m+i,m+j}=0
$
when $j>i$. Then, for $m$ fixed and $i=1,2,\ldots$ we have 
\begin{equation}
\left.
\begin{array}{cccccccccccc}
(\lambda_{m+1}-\lambda_m)b_{m+1,m}& &&&&&& = & \delta_{m+1}^{(1)}\\
(\lambda_{m+2}-\lambda_m)b_{m+2,m}& - &  \delta_{m+1}^{(1)}b_{m+2,m+1}&&&&& = & \delta_{m+2}^{(2)}\\
\vdots & & \vdots &&&&&& \vdots\\
(\lambda_{m+N}-\lambda_m)b_{m+N,m}& - &  \delta_{m+1}^{(1)}b_{m+N,m+1}&-&\cdots & -& \delta_{m+N-1}^{(N-1)}b_{m+N,m+N-1}&= & \delta_{m+N}^{(N)}\\
\vdots & & \vdots &&&&\vdots && \vdots\\
(\lambda_{s}-\lambda_m)b_{s,m}& - &  \delta_{m+1}^{(1)}b_{s,m+1}&-&\cdots & -& \delta_{m+N-1}^{(N-1)}b_{s,m+N-1}&= & \delta_{m+N}^{(N)}b_{s,m+N}\\
\vdots & & \vdots &&&&\vdots && \vdots\\
\end{array}
\right\}
\label{nuevo2}
\end{equation}
Hence, since the first relation in \eqref{nuevo2},
$$
b_{m+1,m}=\frac{ \delta_{m+1}^{(1)}}{\lambda_{m+1}-\lambda_m}
$$
which is \eqref{coeficientes} when $n=m+1$ and $i=m$.

From the second and consecutive relations of \eqref{nuevo2} we obtain \eqref{coeficientes} for $n=m+k,\,i=m,$ and $k=2,3,\ldots$ Assume 
$b_{m+k,m}$ verifying \eqref{coeficientes} for $k=2,3,\ldots, q$. Then since the $(q+1)$-th relation of \eqref{nuevo2} we have
\begin{equation}
\label{nuevo3}
(\lambda_{m+q+1}-\lambda_m)b_{m+q+1,m}=\delta_{m+1}^{(1)}b_{m+q+1,m+1}+\delta_{m+2}^{(2)}b_{m+q+1,m+2}+\cdots+\delta_{m+N}^{(N)}b_{m+q+1,m+N}.
\end{equation}
 
On the other hand, because $m$ can be substituted in \eqref{nuevo2} by any $m+p,\,p\in \mathbb N$, we know that
$$
b_{m+q+1,m+j}=\mathlarger\sum_{E_{m+q+1}^{(m+j)}}\left(\prod_{s=1}^k\frac{\delta^{(i_s)}_{m+j+i_1+\cdots +i_s}}{\lambda_{m+q+1}-\lambda_{m+j+i_1+\cdots +i_{s-1}}}\right),\quad j=1,\ldots , N.
$$
From this and \eqref{nuevo3},
\begin{equation}
b_{m+q+1,m}=\mathlarger\sum_{j=1}^{q+1}\frac{\delta^{(j)}_{m+j}}{\lambda_{m+q+1}-\lambda_{m}}\mathlarger\sum_{E_{m+q+1}^{(m+j)}}\left(\prod_{s=1}^k\frac{\delta^{(i_s)}_{m+j+i_1+\cdots +i_s}}{\lambda_{m+q+1}-\lambda_{m+j+i_1+\cdots +i_{s-1}}}\right),\quad j=1,\ldots , N.
\label{nuevo5}
\end{equation}
Furthermore it is obvious that, if  
$
\sum_{s=1}^ki_s=q+1,
$
then 
$
\sum_{i_s\neq j}i_s=q-j+1
$
for each 
$
j\in \{i_1,\ldots, i_k\}.
$
Thus, it is easy to see
$$
E_{m+q+1}^{(m)}=\bigcup_{j=1}^{q+1} E_{m+q+1}^{(m+j)}.
$$
Consequently, from \eqref{nuevo5} we arrive to
\begin{eqnarray*}
b_{m+q+1,m} &= &\mathlarger\sum_{j=1}^{q+1}\mathlarger\sum_{E_{m+q+1}^{(m+j)}}\frac{\delta^{(j)}_{m+j}}{\lambda_{m+q+1}-\lambda_{m}}\left(\prod_{s=1}^k\frac{\delta^{(i_s)}_{m+j+i_1+\cdots +i_s}}{\lambda_{m+q+1}-\lambda_{m+j+i_1+\cdots +i_{s-1}}}\right)\\
& = & \mathlarger\sum_{E_{m+q+1}^{(m)}}
\left(\prod_{s=1}^k\frac{\delta^{(i_s)}_{m+i_1+\cdots +i_s}}{\lambda_{m+q+1}-\lambda_{m+i_1+\cdots +i_{s-1}}}\right),
\end{eqnarray*}
which is \eqref{coeficientes} for $n=m+q+1$ and $i=m$. \hfill $\square$

For illustrative purposes only, we consider the case of the well-known Hermite polynomials. For them, the differential operator \eqref{1} is
$$
L\equiv a_1(x)\partial_x +a_2(x)\partial_x^2
$$
where
$a_1(x)=2x$ and $a_2(x)=-1$. With some straight computations, from \eqref{A} we obtain
$$
\delta_r^{(k)}=\left\{
\begin{array}{cccc}
2r & \text{if} &k=0\\
0 & \text{if} &k=1\\
-r(r-1) & \text{if} &k=2\\
0 & \text{if} &k\geq 3
\end{array}
\right. ,\qquad r=k,k+1,\ldots 
$$
In particular, the eigenvalues are 
$
\lambda_n=2n,\,n=0,1,\ldots
$
and Theorem \ref{teorema_coeficientes} is applicable because the restriction \eqref{autovectores} is fulfilled.

Due to $\delta_r^{(1)}=0$ and $E_n^{(n-1)}=\{1\}$, in  \eqref{coeficientes} we see $b_{n,n-1}=0 $. Moreover, also when $n-i$ is odd then there exists some $i_s=1$ such that 
$
(i_1,\ldots ,i_k)\in E_n^{(i)}
$
and $b_{n,i}=0 $. Regarding the even coefficients 
$b_{m+2s,m}, m,s,\in \mathbb N$, since 
$E_{m+2s}^{(m)}$
can be substituted by 
$\{\stackrel{(s)}{(2,\ldots,2)}\},$
we have
$$
b_{m+2s,m}=\prod_{j=1}^s\frac{\delta_{m+2j}^{(2)}}{\lambda_{m+2s}-\lambda_{m+2j-2}}=\frac{(-1)^s}{2^{2s}}\frac{(m+2s)!}{m!s!}
$$
or, what is the same,
$$
b_{n,n-2s}=\frac{(-1)^s}{2^{2s}}\frac{n!}{(n-2s)!s!}
$$
Hence we arrive the known monic Hermite polynomials,
$$
P_n(x)=x^n+\sum_{s=1}^{[\frac{n}{2}]}b_{n,n-2s}x^{n-2s}=\frac{n!}{2^n}\sum_{s=0}^{[\frac{n}{2}]}\frac{(-1)^s}{(n-2s)!s!}(2x)^{n-2s}.
$$

Returning to the general case, we are interested in sequences of polynomials verifying a $(p+2)$-term recurrence relation for some $p\in\mathbb N$. Since we are dealing with the coefficients of $P_n(x)$ given in \eqref{111}, in the rest of the section such condition is analyzed in terms of these coefficients. 

For each $n=0,1,\ldots $ we define the infinite matrix
$$
B^{(n)}:=
\begin{pmatrix} 
b_{n,0}&b_{n+1,0}&\cdots&b_{n+k,0}&\cdots \\ 
b_{n,1}&b_{n+1,1}&\cdots&b_{n+k,1}&\cdots\\ 
 \vdots&\vdots &\ddots&\vdots  \\
  b_{n,n-1}&b_{n+1,n-1}&\cdots&b_{n+k,n-1}&\cdots\\ 
1&b_{n+1,n}&\cdots&b_{n+k,n}&\cdots\\ 
0 & 1&   \ddots& \vdots &\ddots &\\
 & 0 & \ddots &b_{n+k,n+k-1}&\cdots\\
& & \ddots & \ddots&\ddots\\
 \end{pmatrix}.
$$
We denote by 
$B^{(n)}_{r,s},\,r,s\in\mathbb{N},$
the main section of order $r\times s$ of $B^{(n)}$. Finally, if $r\geq s$ then
$\widetilde B^{(n)}_{r,s}$ denotes the square matrix of order $s$ formed by the last $s$ rows of $B^{(n)}_{r,s}$. 

\begin{teo}
\label{teorema2}
The sequence of monic polynomials 
$
\{P_n\}
$
verifies the recurrence relation \eqref{otro9} if and only if 
\begin{equation}
B^{(n-p)}_{n-p,p+2}
\left(
\widetilde B^{(n-p)}_{n+2,p+2}
\right)^{-1}
\begin{pmatrix}
b_{n,n-p-1}\\
b_{n,n-p}\\
 \vdots\\
 b_{n,n-1}\\
1
 \end{pmatrix}=
\begin{pmatrix}
0\\
b_{n,0}\\
b_{n,1}\\
\vdots\\
b_{n,n-p-2}
\end{pmatrix}
\label{V}
\end{equation}
for any $n\in \mathbb{N}$ such that $n>p$.
\end{teo}

\noindent
\underline{Proof}.- For each 
$
n=0,1,\ldots, 
$ \eqref{otro9} is verified if and only if 

\begin{equation}
\left.
\begin{array}{rcl}
\displaystyle \sum_{k=n-p}^{n}\alpha_{n,k}b_{k,i}-b_{n,i-1}+b_{n+1,i} & = & 0,\\
 i=0,1,\ldots, n+1\,,
\end{array}
\right\}
\label{4}
\end{equation}
where we understand 
$
b_{n,-1}=0
$.
Then we want to prove that \eqref{V} and \eqref{4} are equivalent.

In the first place, we show that \eqref{V} is a neccesary condition. Using matrix notation, if $n\geq p$ then \eqref{4} is rewritten as 
\begin{equation}
B^{(n-p)}_{n+2,p+2}
\begin{pmatrix}
 \alpha_{n,n-p}\\
 \alpha_{n,n-p+1}\\
 \vdots\\
 \alpha_{n,n}\\
1
 \end{pmatrix}=
\begin{pmatrix}
0\\
b_{n,0}\\
b_{n,1}\\
\vdots\\
b_{n,n-1}\\
1
\end{pmatrix}
\label{I}
\end{equation}

If
$
(\alpha_{n,n-p}, \ldots, \alpha_{n,n},1)
$
verifies \eqref{I} then 
\begin{equation}
\widetilde B^{(n-p)}_{n+2,p+2}
\begin{pmatrix}
 \alpha_{n,n-p}\\
 \alpha_{n,n-p+1}\\
 \vdots\\
 \alpha_{n,n}\\
1
 \end{pmatrix}=
\begin{pmatrix}
b_{n,n-p-1}\\
b_{n,n-p}\\
\vdots\\
b_{n,n-1}\\
1
\end{pmatrix}
\label{II}
\end{equation}
and 
\begin{equation}
B^{(n-p)}_{n-p,p+2}
\begin{pmatrix}
 \alpha_{n,n-p}\\
 \alpha_{n,n-p+1}\\
 \vdots\\
 \alpha_{n,n}\\
1
 \end{pmatrix}=
\begin{pmatrix}
0\\
b_{n,0}\\
b_{n,1}\\
\vdots\\
b_{n,n-p-2}
\end{pmatrix}\,.
\label{III}
\end{equation}
(Note that $\widetilde B^{(n-p)}_{n+2,p+2}$ is a triangular matrix.)
In fact, \eqref{I} is equivalent to \eqref{II}-\eqref{III}. Moreover,  \eqref{II}-\eqref{III} can be expressed as \eqref{V}, as we wanted to prove.

On the other hand, for $n\leq p$, the coefficients $\alpha_{n,i},\,i=0,1,\ldots, n,$ of \eqref{4} alway are defined by
\begin{equation}
B^{(0)}_{n+2,n+2}\begin{pmatrix}
 \alpha_{n,0}\\
 \alpha_{n,1}\\
 \vdots\\
 \alpha_{n,n}\\
1
 \end{pmatrix}=
\begin{pmatrix}
0\\
b_{n,0}\\
b_{n,1}\\
\vdots\\
b_{n,n-1}\\
1
\end{pmatrix},
\label{IV}
\end{equation}
because
\begin{equation}
B^{(0)}_{n+2,n+2}=
\begin{pmatrix}
 1 & b_{1,0} & \cdots &\cdots & b_{n+1,0}\\
  & 1 & b_{2,1} & \ddots & \vdots \\
& &\ddots &\ddots &  \vdots \\
&&&\ddots &\vdots \\
&& &  &1
 \end{pmatrix}
\label{triangular}
\end{equation}
is a triangular regular matrix. Then,  \eqref{V} is also a sufficient condition to have \eqref{otro9} for each $n\in \mathbb{N}$.

\hfill $\square$

\begin{coro}
\label{corolario} 
If the sequence 
$
\{P_n\}
$
verifies \eqref{otro9} then we have

\begin{equation}
\label{19*}
\begin{pmatrix}
 \alpha_{n,0}\\
 \alpha_{n,1}\\
 \vdots\\
 \alpha_{n,n}\\
1
 \end{pmatrix}=\left(B^{(0)}_{n+2,n+2}\right)^{-1}
\begin{pmatrix}
0\\
b_{n,0}\\
b_{n,1}\\
\vdots\\
b_{n,n-1}\\
1
\end{pmatrix},\quad \text{if } n\leq p,
\end{equation}
and 
\begin{equation}
\label{19**}
\begin{pmatrix}
 \alpha_{n,n-p}\\
 \alpha_{n,n-p+1}\\
 \vdots\\
 \alpha_{n,n}\\
1
 \end{pmatrix}=\left(\widetilde B^{(n-p)}_{n+2,p+2}\right)^{-1}
\begin{pmatrix}
b_{n,n-p-1}\\
\vdots\\
b_{n,n-1}\\
1
\end{pmatrix}, \quad \text{if } n> p.
\end{equation}
\end{coro}

\noindent
\underline{Proof}.- This is an immediate consequence of \eqref{II} and \eqref{IV}.

In the proof of Theorem \ref{teorema2} we see that 
$
\{P_n\}
$
verify \eqref{otro9} if and only if \eqref{I} takes place. Besides, \eqref{I} is equivalent to 
\begin{equation}
B^{(n-p)}_{n+1,p+1}
\begin{pmatrix}
 \alpha_{n,n-p}\\
 \alpha_{n,n-p+1}\\
 \vdots\\
 \alpha_{n,n}
 \end{pmatrix}
=
\begin{pmatrix}
0\\
b_{n,0}\\
b_{n,1}\\
\vdots\\
b_{n,n-1}
\end{pmatrix}-
\begin{pmatrix}
b_{n+1,0}\\
b_{n+1,1}\\
\vdots\\
b_{n+1,n}
\end{pmatrix}
\label{I*}
\end{equation}
for each $n\geq p$, because 
$
B^{(n-p)}_{n+2,p+1}$
is obtained by eliminating the last column from 
$
B^{(n-p)}_{n+2,p+2}$. This means that 
$
\left(
x_1,x_2,\ldots ,x_{p+1}
\right)=
\left(
 \alpha_{n,n-p},  \alpha_{n,n-p+1}, \ldots ,\alpha_{n,n}
\right)
$
is a solution of the following overdetermined linear system,
\begin{equation}
\label{sobredeterminado}
\left.
\begin{array}{ccl}
b_{n-p,i}x_1+b_{n-p+1,i}x_2+\cdots + b_{n,i}x_{p+1}& = &b_{n,i-1}-b_{n+1,i}\\\\
&i=&0,1,\ldots ,n, 
\end{array}
\right\}
\end{equation}
what is the same that \eqref{4} (we recall that 
$
b_{i,j}=0
$
when
$i<j$).
From this and Corollary \ref{corolario} we see that the solution is unique when it exists, and we can explicitly obtain these coefficients $\{\alpha_{n,n-k}\}$ of the recurrence relation (see \eqref{19*}). Moreover, as a consequence of \eqref{4} and Remark \ref{remark 1}, they are rational functions of $n$. Therefore, proving the existence of the searched recurrence relation is reduced to verifying \eqref{19**} or, equivalently, checking \eqref{sobredeterminado} for 
$n>p$ and
$i=0,1,\ldots, n-p-1.$

On the other hand, for each fixed $n\in \mathbb{N}$, \eqref{matriz} is a recurrence relation with the $N+1$ terms 
$$b_{n,m},b_{n,m+1}\ldots, b_{n,m+N},\,m=0,1,\ldots$$
Therefore, for $n>N$ it is possible to state the coefficients $b_{n,i}$ of each polynomial $P_n(x)$ in terms of the first $N$ coefficients (or either in terms of the last ones $N$).  This suggests the existence of a $(N+1)$-term recurrence relation for the sequence of eigenpolynomials $\{P_n\}$. However, note that some of these terms can be null and then have a recurrence relation with a smaller number of terms. 

As a summary of the above reflections, we establish the following

\begin{conje}
\label{conjetura}
Given the differential operator \eqref{1} and the sequence
$\{P_n\}$
of eigenpolynomials, there exist $p\leq N-1$ and a sequences of coefficients
$
\{\alpha_{n,k}\}
$ such that \eqref{otro9} is verified.
\end{conje}

\section{The inverse problem\label{seccion4}} 

Now we assume that the sequences $\{\lambda_n\}$ and $\{P_n\}$ are known. In this section our goal is to find 
polynomials 
$a_i(x)$
and 
$N\in\mathbb N$
such that 
$L$ 
is defined as in \eqref{1}, and \eqref{11} takes place.

In the first place, for each 
$n\in \mathbb{N}$
we  will define
$
\delta_n^{(0)}, \,\delta_n^{(1)},\ldots, \delta_n^{(n)}
$
verifying \eqref{7}. Take 
\begin{equation}
\label{auto}
\delta_n^{(0)}=\lambda_n,\quad n=0,1,\ldots 
\end{equation}
(we recall that we assume 
$\lambda_0=0$).
Then \eqref{7} holds for $m=n$.  We define the rest of the values as
\begin{equation}
\label{1*}
\delta_n^{(k)}=(-1)^k
\left|
\begin{array}{cccccc}
 (\lambda_{n-k}-\lambda_{n-k+1}) b_{n-k+1,n-k}& (\lambda_{n-k}-\lambda_{n-k+2}) b_{n-k+2,n-k}&\cdots&\cdots  &(\lambda_{n-k}-\lambda_{n}) b_{n,n-k}\\ 
1& b_{n-k+2,n-k+1}&\cdots &\cdots& b_{n,n-k+1}\\ 
0 & 1&   \cdots& \cdots &b_{n,n-k+2} &\\
\vdots  & 0  & \ddots &&\vdots \\
\vdots  &\vdots & \ddots &\ddots &\vdots \\
0 & 0 &\cdots & 1& b_{n,n-1}
 \end{array}
\right|
\end{equation}
for each $n\in \mathbb{N}$ and $k=1,2,\ldots,n.$ The following result guarantees that the sequence
$
\{
\delta_n^{(k)}
\}
$
can be used to construct $L$. 

\begin{teo}
\label{teorema3}
For all $n\in \mathbb{N}$ and $k=1,2,\ldots,n$, 
$
\{
\delta_n^{(k)}
\}
$
verifies 
\begin{equation}
\label{otra10}
\sum_{k=0}^{n-m}\delta_{m+k}^{(k)}b_{n,m+k}=\lambda_nb_{n,m}\,,\quad m=0,1,\ldots, n\,.
\end{equation}

\end{teo}

\noindent
\underline{Proof}.- We proceed by induction on $n$. For $n=1$ we have 
$
\delta_1^{(0)}=\lambda_1$
in \eqref{auto} and $\delta_1^{(1)}=\lambda_1b_{1,0}
$
in  \eqref{1*}, which is \eqref{otra10} for $n=1$ and, respectively, $m=1,\,0$. 

Let 
$
n\in\mathbb{N},\,n>1,$
be and assume that \eqref{otra10} takes place for $m=0,1,\ldots,n$. Substituting $n$ by $n+1$, \eqref{otra10} becomes 
\begin{equation}
\label{nuevo*}
\sum_{k=0}^{n-m+1}\delta_{m+k}^{(k)}b_{n+1,m+k}=\lambda_{n+1}b_{n+1,m}\,,\quad m=0,1,\ldots, n+1\,,
\end{equation}
which we want to prove. 

For 
$m=n+1$ \eqref{nuevo*} is 
$$
\delta_{n+1}^{(0)}=\lambda_{n+1},
$$
that we know since \eqref{auto}.

For $m=0,\ldots,n$, taking $s=n-m$ we write \eqref{nuevo*} as
\begin{equation}
\label{nuevo**}
\sum_{k=0}^{s+1}\delta_{n-s+k}^{(k)}b_{n+1,n-s+k}=\lambda_{n+1}b_{n+1,n-s}\,,\quad s=0,1,\ldots, n\,.
\end{equation}
in \eqref{nuevo*}.
Developing the determinant
\begin{equation}
\label{2*}
\delta_{n+1}^{(k)}=(-1)^k
\left|
\begin{array}{cccccc}
 (\lambda_{n-k+1}-\lambda_{n-k+2}) b_{n-k+2,n-k+1}&\cdots&\cdots  &(\lambda_{n-k+1}-\lambda_{n+1}) b_{n+1,n-k+1}\\ 
1& \cdots &\cdots& b_{n+1,n-k+2}\\ 
0 & &   \cdots& b_{n+1,n-k+3} &\\
\vdots  & & \ddots &\vdots \\
0& &\cdots &  b_{n+1,n}
 \end{array}
\right|
\end{equation}
by its last column we arrive to 
\begin{equation}
\label{otro25}
\delta_{n+1}^{(s+1)}= (\lambda_{n+1}-\lambda_{n-s})b_{n+1,n-s}-\sum_{k=1}^{s}\delta_{n+k-s}^{(k)}b_{n+1,n+k-s},
\end{equation}
which is \eqref{nuevo**}.
\hfill $\square$

\begin{coro}
In the above conditions, if there exists 
$N\in \mathbb{N}$
such that 
\begin{equation}
\label{nuevo***}
\delta_n^{(k)}=0\, \text{   for each }n>N,\,n \in \mathbb N, \text{ and } k=N+1,\ldots , n,
\end{equation}
then there exists a differential operator $L$ as in \eqref{1} such that \eqref{11} takes place. 
\end{coro}

\noindent
\underline{Proof}.- Since \eqref{nuevo***}, in \eqref{otra10} we have \eqref{7}. From the sequence 
$
\{\delta_n^{(k)}:n\in \mathbb N, \,k=0,1,\ldots, n\}
$
we define the sequence of polynomials 
$
\{a_m(x)\},\,m\in \mathbb N,
$
whose coefficients are given in \eqref{B} (see also \eqref{1111}). Then the differential operator $L$ is defined as in \eqref{1} from 
$
a_1(x),\ldots, a_N(x).
$
In fact, following the proof of Theorem \ref{lema1} we see that another way to write \eqref{11} is \eqref{22}, from which \eqref{7} is obtained by comparing coefficients.  
\hfill $\square$

In the above corollary we have defined the polynomials 
$a_n(x)$
for any 
$n\in \mathbb N$,
despite only those that satisfy 
$\deg(a_n)\leq N
$
are used in the construction of $L$. We underline that condition \eqref{nuevo***} implies 
$$
n!a_{n,n-k}=\sum_{i=k}^N {n \choose i}(-1)^{n-i}\delta_i^{(k)}\quad \text{ for } k=0,1,\ldots, n\,,
$$
and, in particular, 
$
a_{n,n-k}=0
$
when
$
k=N+1,\ldots, n.
$
This is, 
$$
a_n(x)=\sum_{k=0}^Na_{n,n-k}x^{n-k}=x^{n-N}\sum_{i=n-N}^na_{n,i}x^{i-n+N}.
$$ 
In other words, 
$
\displaystyle\frac{a_n(x)}{x^{n-N}}
$
is a polynomial such that
$$
\deg\left(
\frac{a_n(x)}{x^{n-N}}
\right)\leq N, \quad \forall n>N.
$$

Note that we are starting from any sequence of monic polynomials 
$
\{P_n\}
$, 
which may not verify a finite recurrence relation \eqref{otro9}. However, if \eqref{nuevo***} takes place and Conjecture \ref{conjetura} is true, then necessarily \eqref{otro9} is fulfilled for $p=N-1$.

\end{document}